\RequirePackage{ifpdf}
\ifpdf 
\documentclass[pdftex]{sigma}
\else
\documentclass{sigma}
\fi

\numberwithin{equation}{section}

\newcommand\idd{\mathop {\fam 0 id}\nolimits}
\newcommand\Curr{\mathop {\fam 0 Cur}\nolimits}
\newcommand\Ann{\mathop {\fam 0 Ann}\nolimits}
\newcommand\eval{\mathop {\fam 0 ev}\nolimits}
\newcommand\Cend{\mathop {\fam 0 Cend}\nolimits}
\newcommand\Chom{\mathop {\fam 0 Chom}\nolimits}

\newcommand\codim{\mathop {\fam 0 codim}\nolimits}
\newcommand\Coeff{\mathop {\fam 0 Coef\/f}\nolimits}
\newcommand\Derr{\mathop {\fam 0 Der}\nolimits}
\newcommand\End{\mathop {\fam 0 End}\nolimits}

\newcommand{\oo}[1]{\mathrel{{\circ }_{#1}} }

\begin{document}

\allowdisplaybreaks

\renewcommand{\thefootnote}{$\star$}

\renewcommand{\PaperNumber}{014}

\FirstPageHeading

\ShortArticleName{Simple Finite Jordan Pseudoalgebras}

\ArticleName{Simple Finite Jordan Pseudoalgebras\footnote{This paper is a
contribution to the Special Issue on Kac--Moody Algebras and Applications. The
full collection is available at
\href{http://www.emis.de/journals/SIGMA/Kac-Moody_algebras.html}{http://www.emis.de/journals/SIGMA/Kac-Moody{\_}algebras.html}}}

\Author{Pavel KOLESNIKOV}

\AuthorNameForHeading{P.~Kolesnikov}

\Address{Sobolev Institute of Mathematics, 4 Acad. Koptyug Ave.,
630090 Novosibirsk, Russia}
\Email{\href{mailto:pavelsk@math.nsc.ru}{pavelsk@math.nsc.ru}}

\ArticleDates{Received September 12, 2008, in f\/inal form January 10,
2009; Published online January 30, 2009}

\Abstract{We consider the structure of
Jordan $H$-pseudoalgebras which are linearly f\/initely
generated over a Hopf algebra~$H$.
There are two cases under consideration:
$H=U(\mathfrak h)$
and
$H=U(\mathfrak h)\# \mathbb C[\Gamma ]$,
where $\mathfrak h$ is
a f\/inite-dimensional Lie algebra over $\mathbb C$,
$\Gamma $ is an arbitrary group acting on
$U(\mathfrak h)$ by automorphisms.
We construct
an analogue of the Tits--Kantor--Koecher
construction for f\/inite Jordan pseudoalgebras
and describe all simple ones.}


\Classification{17C50; 17B60; 16W30} 


\section{Introduction}\label{sec1}

The notion of pseudoalgebra
appeared as a natural generalization
of the notion of conformal algebra.
The last one provides a formal language describing
algebraic structures
underlying the singular part of the operator product
expansion (OPE)
in conformal f\/ield theory. Roughly speaking,
the OPE of two local chiral f\/ields is a formal distribution in
two variables presented as
$\sum\limits_{n=-\infty}^{N-1}  c_n (z) (w-z)^{-n-1} $
\cite{BPZ}.
The coef\/f\/icients $c_n$, $n\in \mathbb Z$,
of this distribution are considered
as new ``products'' on the space of f\/ields.
The algebraic system
obtained is called a vertex algebra.
Its formal axiomatic description was stated
in~\cite{Bor} (see also~\cite{K1}).
The ``singular part'' of a vertex algebra, i.e.,
the structure def\/ined by only those operations
with non-negative~$n$,  is a~(Lie) conformal algebra
\cite{K1}.

Another approach to the theory of vertex algebras
gives rise to the notion of a pseudotensor category~\cite{BD} (which is similar to the multicategory
of~\cite{La}).
Given a Hopf algebra $H$, one may def\/ine
the pseudotensor category $\mathcal M^*(H)$
\cite{BDK} (objects of this category are left $H$-modules).
An algebra in this category is called
an $H$-pseudoalgebra. A pseudoalgebra is said to be f\/inite
if it is a f\/initely generated $H$-module.

In particular, for the
one-dimensional Hopf algebra
$H=\Bbbk$, $\Bbbk $ is a f\/ield, an $H$-pseudoalgebra
is just an ordinary algebra over the f\/ield~$\Bbbk $.
For $H=\Bbbk[D]$, an $H$-pseudoalgebra is
exactly the same as conformal algebra.
In a more general case $H=\Bbbk[D_1,\dots, D_n]$,
$n\ge 2$, the notion of an $H$-pseudoalgebra
is closely related with Hamiltonian
formalism in the theory of non-linear evolution
equations \cite{BDK}. For an arbitrary Hopf algebra $H$,
an $H$-pseudoalgebra def\/ines a functor
from the category of $H$-bimodule
associative commutative algebras
to the category of
$H$-module algebras
(also called $H$-dif\/ferential algebras).

An arbitrary conformal algebra $C$ can be canonically embedded
in a ``universal'' way into
the space of formal power series $A[[z,z^{-1}]]$ over an
appropriate ordinary algebra $A$ \cite{K1, Ro1}. This algebra
$A=\Coeff C$ is called the coef\/f\/icient algebra of $C$. A conformal algebra
is said to be associative (Lie, Jordan, etc.) if so is its
coef\/f\/icient algebra. For pseudoalgebras, a~construction called
annihilation algebra \cite{BDK} works instead of coef\/f\/icient algebra.
However, the notion of a~pseudotensor
category provides a direct way to the def\/inition of what is a variety of
pseudoalgebras~\cite{Kol}.

In the paper
\cite{DK},
the complete description of simple
f\/inite Lie conformal algebras over $\Bbbk=\mathbb C$
was obtained. Apart from current conformal algebras,
the only example of a simple f\/inite Lie conformal
algebra is the Virasoro conformal algebra.
In the associative case,
there are no exceptional
examples: A simple f\/inite
associative conformal algebra
is isomorphic to the current algebra over
$M_n(\mathbb C)$, $n\ge 1$.
It was shown in
\cite{Z2000} that
there are no exceptional examples of
simple f\/inite Jordan conformal algebras.

In \cite{BDK}, the structure theory of
f\/inite Lie pseudoalgebras was developed.
The classif\/ication theorem \cite[Theorem~13.2]{BDK}
states that there
exist simple f\/inite Lie pseudoalgebras
which are not isomorphic to current pseudoalgebras
over ordinary simple f\/inite-dimensional Lie algebras.
This is not the case for associative
pseudoalgebras. In this paper, we show
the same for f\/inite simple Jordan pseudoalgebras
(Theorem~\ref{thm 5.1}): There are no examples of such
pseudoalgebras except for
current algebras (if $H=U(\mathfrak h)$)
or
transitive direct sums of current algebras
(if $H=U(\mathfrak h)\#
\mathbb C [\Gamma]$). The main tool of the proof is an analogue of
the well known Tits--Kantor--Koecher (TKK) construction for Jordan
algebras. This result generalizes the classif\/ication of f\/inite
Jordan conformal algebras
\cite{Z2000}
to ``multi-dimensional'' case.

It was shown in \cite{KR} that the structure theory
 of Jordan conformal superalgebras
is richer.  The classif\/ication of simple f\/inite Jordan superalgebras
based on the structure theory of f\/inite conformal Lie superalgebras
\cite{FK} includes
one series and two exceptional algebras \cite[Theorem~3.9]{KR}.

In our proof, we will not use annihilation algebras
directly, the TKK construction will be built on
the level of pseudoalgebras.

The paper is organized as follows.
Section~\ref{sec2} contains the
basics of Hopf algebras and pseudoalgebras theory,
and some notations that will be used later.
In Section~\ref{sec4}, we introduce an analogue of the
Tits--Kantor--Koecher construction
(TKK) for f\/inite Jordan pseudoalgebras.
To complete the classif\/ication of f\/inite
simple Jordan pseudoalgebras,
we need some technical results proved in
Section~\ref{sec5}.
 The main case under consideration is
$H=U(\mathfrak h)$,
where
$\mathfrak h$
is a~f\/inite-dimensional Lie algebra over~$\mathbb C$.
Another case is the smash-product $U(\mathfrak h)\# \mathbb C[\Gamma]$,
where $\Gamma $ is an arbitrary group.
These cases describe all
cocommutative Hopf algebras over $\mathbb C$ with
f\/inite-dimensional spaces of primitive elements
(see, e.g., \cite{Sw}).

\section{Preliminaries on Hopf algebras and pseudoalgebras}\label{sec2}

\subsection{Hopf algebras}\label{subsec2.1}

In this section, we state some
notations that will be used later.

An associative algebra
 $H$  with a unit (over  a f\/ield $\Bbbk$)
endowed with coassociative coproduct
$\Delta : H \to H\otimes H$
and counit
$\varepsilon : H\to \Bbbk$
is called a bialgebra. Recall that
both $\Delta$
and
$\varepsilon$ are homomorphisms of algebras and
\[
(\idd \otimes \Delta)\Delta(h) = (\Delta \otimes \idd)\Delta(h),
\qquad
(\varepsilon\otimes\idd)\Delta(h)
 = (\idd \otimes \varepsilon)\Delta(h) = h.
\]
To simplify the notation, we will  use the following one which
is due to Sweedler~\cite{Sw}:
$\Delta^{[1]}(h):=h$,
$\Delta^{[2]}(h) := \Delta (h)=
\sum\limits_{(h)} h_{(1)} \otimes h_{(2)}$,
$\Delta^{[n]}(h) := (\idd\otimes \Delta^{[n-1]})
= \sum\limits_{(h)}h_{(1)} \otimes \dots \otimes h_{(n)}$.
Further, we will omit the symbol $\sum\limits_{(h)} $ by writing
$\Delta(h)
  = h_{(1)}\otimes h_{(2)}$,
$\Delta^{[n]}(h)= h_{(1)} \otimes \dots \otimes h_{(n)}$,
  etc.

Given a bialgebra $H$, an antihomomorphism
$S : H\to H$
is called an antipode, if it satisf\/ies
\[
S(h_{(1)})h_{(2)} =
 \varepsilon (h)= h_{(1)}S(h_{(2)}).
\]
A bialgebra with an antipode is called a Hopf algebra.

There exists a natural structure of (right)
$H$-module on the $n$th
tensor power of $H$ (denoted by $H^{\otimes n}$):
\begin{equation}
(f_1\otimes \dots \otimes f_n)\cdot
h = f_1 h_{(1)}\otimes\dots \otimes f_n h_{(n)},
\qquad
f_i, h\in H.
                                               \label{eq2.2}
\end{equation}

In this paper, we substantially consider
cocommutative Hopf algebras,
i.e., such that
$h_{(1)}\otimes h_{(2)}= h_{(2)}\otimes h_{(1)}$
for all $h\in H$. The antipode $S$ of a cocommutative
Hopf algebra is involutive, i.e., $S^2=\idd $.

For example,
the universal enveloping algebra $U(\mathfrak g)$
of a Lie algebra $\mathfrak g$ over a f\/ield of zero
characteristic is a cocommutative Hopf algebra.
Another series of examples is provided by
the group algebra $\Bbbk[\Gamma ]$ of an arbitrary group $\Gamma $
and by the general construction of smash-product.
Namely, suppose $H$ is a Hopf algebra, and
let a group $\Gamma $ acts on
 $H$ by algebra automorphisms.
Then one may def\/ine the following new product on
$H\otimes \Bbbk[\Gamma ]$:
\[
(h_1\otimes g_1)\cdot (h_2\otimes g_2)
  = h_1h_2^{g_1}\otimes g_1g_2.
\]
The algebra obtained is denoted by $H\# \Bbbk[\Gamma ]$.
Together with usual coproduct and antipode def\/ined as
 on $H\otimes \Bbbk[\Gamma ]$,
$H\# \Bbbk[\Gamma ]$ is a Hopf algebra (the smash product of
$H$ and $\Bbbk [G]$).
If $H$ is cocommutative, then so is $H\# \Bbbk[\Gamma]$.
Moreover, if  $\Bbbk $ is an algebraically closed f\/ield of
zero characteristic, then every cocommutative Hopf algebra $H$
over $\Bbbk $
is isomorphic to the smash product
$U(\mathfrak g)\# \Bbbk[\Gamma ]$ for appropriate $\mathfrak g$ and $\Gamma $
\cite{Sw2}.

\begin{lemma}[\cite{BDK}]\label{lem 2.1}
Let $H$ be a cocommutative Hopf algebra, and
let $\{h_i\mid i\in I\}$ be a linear basis of~$H$.
Then every element
$F\in H^{\otimes n+1}$, $n\ge 1$,
 can be uniquely presented as
\begin{equation}
F=\sum\limits_{i_1,\dots, i_n}
 (h_{i_1}\otimes \dots \otimes h_{i_n} \otimes 1)\cdot g_{i_1,\dots, i_n} ,
 \qquad g_{i_1,\dots, i_n} \in H.
                                         \label{eq2.4}
\end{equation}
\end{lemma}

In other words, the set
$\{h_{i_1}\otimes \dots \otimes h_{i_n}\otimes 1 \mid i_1,\dots , i_n\in I\}$
is an $H$-basis of the $H$-module
$H^{\otimes n+1}$~\eqref{eq2.2}.

To f\/ind the presentation \eqref{eq2.4}, one may use formal Fourier
transformation $\mathcal F$ and its inver\-se~$\mathcal F^{-1}$~\cite{BDK}:
\begin{gather*}
 \mathcal F, \mathcal F^{-1} : \ H^{\otimes n+1} \to H^{\otimes n+1},\nonumber \\
 \mathcal F: \ h_{1}\otimes \dots \otimes h_{n} \otimes f
 \mapsto h_{1}f_{(1)}\otimes \dots \otimes h_{n}f_{(n)} \otimes f_{(n+1)},
\\
 \mathcal F^{-1}: \ h_{1}\otimes \dots \otimes h_{n} \otimes f
 \mapsto h_{1}S(f_{(1)})\otimes \dots \otimes h_{n}S(f_{(n)}) \otimes f_{(n+1)}.
\end{gather*}
We will use a ``left'' analogue of the Fourier transformation
\[
\mathcal F': \  h\otimes f_1\otimes \dots \otimes f_n \mapsto
  h_{(1)} \otimes h_{(2)} f_1 \otimes \dots \otimes  h_{(n+1)}f_n,
\]
which is also invertible.

\subsection{Dual algebras}\label{subsec2.2}

Suppose
$H$ is a cocommutative Hopf algebra,
and let
$X=H^*$ be its dual algebra (i.e., the product on $X$
is dual to the coproduct on $H$).
Let us f\/ix a linear basis  $\{h_i \mid i\in I\}$ of $H$
and denote by
$\{x_i \mid i\in I\}\subset X$
 the set of dual functionals:
$\langle x_i, h_j \rangle = \delta_{ij}$,
$i,j\in I$.

An arbitrary element $x\in X$
can be presented as an inf\/inite series in $x_i$:
\[
x = \sum\limits_{i\in I}\langle x, h_i\rangle x_i.
\]
The algebra $X$ is a left and right module over $H$ with respect to
the actions given
by
\begin{equation}
\langle xh, g\rangle =\langle x, gS(h)\rangle ,
\qquad
\langle hx, g\rangle = \langle x, S(h)g \rangle ,
\qquad
x\in X,\ h,g\in H.
                                     \label{eq2.5}
\end{equation}
The actions \eqref{eq2.5} turn $X$ into a dif\/ferential $H$-bimodule,
i.e.,
$(xy)h= (xh_{(1)})(y h_{(2)})$,
$h(xy)= (h_{(1)}x)(h_{(2)}y)$.

The operation
$\Delta_X : X \to \overline{X\otimes X}:= (H\otimes H)^*$
dual to the product on $H$ is somewhat similar to a coproduct.
From the combinatorial point of view,
$\overline{X\otimes X}$
can be considered as the linear space
of all inf\/inite series
$\sum\limits_{i,j\in I}
 \alpha_{ij}x_i\otimes x_j$,
$\alpha_{ij}\in \Bbbk$.

In order to unify notations, we will use
$x_{(1)}\otimes x_{(2)}$ for $\Delta_X(x)$, $x\in X$.
In particular, the analogues of Fourier transforms
\[
\mathcal F: \ x\otimes y \mapsto xy_{(1)}\otimes y_{(2)},
\qquad
\mathcal F^{-1}: x\otimes y \mapsto xS^*(y_{(1)})\otimes y_{(2)},
\qquad x,y\in X,
\]
act from $X\otimes X$ to $\overline{X\otimes X}$.

\begin{definition}\label{defn 2.2}
Let $V$
be a linear space.
A linear map
 $\pi : X\otimes X \to V$
is said to be {\it local},
if
\[
\pi(x_i\otimes x_j) = 0 \quad
  \hbox{\rm for almost all of the pairs}
\quad (i,j) \in  I^2.
\]
\end{definition}

A local map
$\pi : X\otimes X \to V$
can be naturally continued to
$\bar \pi : \overline{X\otimes X} \to V$. The  map
$\bar\pi $ is continuous with respect to the topology on
$\overline{X\otimes X}$ def\/ined by the following family of
basic neighborhoods of zero:
\[
U^{\perp}=\{\xi \in \overline{X\otimes X} \mid \langle \xi, U\rangle =0\},
\qquad
U\subseteq H\otimes H,\qquad \dim U<\infty
\]
(we assume $V$ is endowed with discrete topology).
Conversely, given a continuous linear map $\overline{X\otimes X}\to V$, its restriction
to $X\otimes X$ is local.

For example, let us f\/ix $h_1,h_2\in H$
and consider the map
$ x\otimes y \mapsto
 \langle x, h_1 \rangle
 \langle y, h_2 \rangle $.
It is clear that this map is local.
Obviously, every local map
$\pi : X\otimes X \to V$
is actually the ``evaluation'' map
\begin{equation}
 \pi(x\otimes y) =
  \eval_a(x,y):=
(\langle x, \cdot \rangle
 \otimes
\langle y, \cdot \rangle
\otimes \idd_V)(a)
                                             \label{eq2.8}
\end{equation}
for an appropriate $a\in H\otimes H\otimes V$.

\begin{lemma}\label{lem 2.3}
Suppose $\pi :X\otimes X \to V$ is a local map, and denote
$\pi\mathcal F=\bar \pi\circ \mathcal F$,
$\pi\mathcal F^{-1} =\bar \pi\circ \mathcal F^{-1}$.
Then both $\pi\mathcal F$ and $\pi\mathcal F^{-1}$ are local,
$\pi\mathcal F = 0$ implies $\pi=0$, and $\pi\mathcal F^{-1} = 0$
implies $\pi =0$.
\end{lemma}

\begin{proof}
Formally speaking, we can not use $\mathcal F^{-1}$ as an inverse of
$\mathcal F$ since both $\mathcal F$ and $\mathcal F^{-1}$ are not def\/ined on
the entire space $\overline{X\otimes X}$.
But it is straightforward to check (see also \cite{BDK}) that
\[
\Delta_X(x) = x_{(1)}\otimes x_{(2)}
 = \sum\limits_{i\in I} xS(h_i)\otimes x_i
=\sum\limits_{i\in I} x_i\otimes S(h_i)x, \quad x\in X,
\]
so
\begin{gather*}
\pi\mathcal F(x_i\otimes x_j)
 =  \bar \pi
\left(
 \sum\limits_{k\in I}x_i(x_jS(h_k))\otimes x_k
\right)
 = \bar \pi
\left (
 \sum\limits_{k,l\in I}\langle x_i(x_jS(h_k)), h_l\rangle
  x_l\otimes x_k
\right)
    \\
\phantom{\pi\mathcal F(x_i\otimes x_j)}{}  = \sum\limits_{k,l\in I}
 \langle x_i, h_{l(1)}\rangle
 \langle  x_j, h_{l(2)}h_k \rangle
 \pi (x_l\otimes x_k).
                                \nonumber
\end{gather*}
This is easy to deduce that if
$\pi = \eval_a$ as in \eqref{eq2.8}
then
\[
\pi \mathcal F(x\otimes y) = \pi'(x\otimes y),
\qquad
\pi'= \eval_{a'},
\qquad
a'=(\mathcal F'\otimes \idd_V)(a)\in
 H\otimes H\otimes V.
\]
Since $\mathcal F': H\otimes H\to H\otimes H$
is invertible, $a=0$ if\/f $a'=0$.
Hence, $\pi'=\pi\mathcal F=0$ implies $\pi=0$.

For $\pi\mathcal F^{-1}$ the proof is completely analogous.
\end{proof}

\subsection{Pseudoalgebras}\label{sec3}

In the exposition of the notion of pseudoalgebra we will preferably follow
\cite{BDK}.

Hereinafter, $H$ is a
cocommutative Hopf algebra,
e.g.,  $H=U(\mathfrak g)$
or $H= U(\mathfrak g)\# \Bbbk [\Gamma ]$.

\begin{definition}[\cite{BDK}]\label{defn 3.1}
Let $P$ be a left $H$-module.
A {\it pseudoproduct\/} is an $H$-bilinear map
\[
*: \ P\otimes P \to (H\otimes H)\otimes_H P.
\]

An $H$-module $P$ endowed with a pseudoproduct~$*$
is called
a {\it pseudoalgebra\/} over $H$ (or $H$-{\it pseudoalgebra}).
If $P$ is a f\/initely generated $H$-module, then
$P$ is said to be {\it finite\/} pseudoalgebra.
\end{definition}

For every $n,m\ge 1$, an $H$-bilinear map
$*: P\otimes P\to (H\otimes H)\otimes _H P$
can be naturally expanded to a map from
$(H^{\otimes n}\otimes_H P)\otimes (H^{\otimes m}\otimes_H P)$
to
$(H^{\otimes n+m}\otimes_H P)$:
\begin{equation}
(F\otimes_H a)*(G\otimes _H b) = ((F\otimes G)\otimes_H 1)
((\Delta^{[n]}\otimes \Delta^{[m]})\otimes_H\idd_P)(a*b),
                                               \label{eq3.2}
\end{equation}
where
$F\in H^{\otimes n}$,
$G\in H^{\otimes m}$,
$a,b\in P$.

This operation allows  to consider long terms in~$P$ with
 respect to~$*$.

One of the main features of a cocommutative bialgebra $H$
is that symmetric groups $S_n$ act by
$H$-module automorphisms
on $H^{\otimes n}$ with respect to
(\ref{eq2.2}).
The action of $\sigma \in S_n$ is def\/ined by
\[
\sigma(h_1\otimes \dots \otimes h_n) =
h_{1\sigma^{-1}} \otimes \dots \otimes h_{n\sigma^{-1}}.
\]

Let us write down the obvious rules matching the action of
$S_n$ with the ``expanded'' pseudoproduct \eqref{eq3.2}.
For every $A\in H^{\otimes n}\otimes_H P$,
$B\in H^{\otimes m}\otimes_H P$,
$\tau\in  S_n$, $\sigma \in S_m$ we have
\begin{equation}
((\tau \otimes_H \idd_P)(A))* B =
(\bar\tau\otimes_H\idd_P)(A*B),
                              \label{eq:3.3}
\end{equation}
where
 $\bar\tau\in S_{n+m}$,
$k\bar\tau =k\tau $ for $k=1,\dots ,n$,
$k\bar\tau =k$ for $k=n+1,\dots ,n+m$;
\begin{equation}
A* ((\sigma  \otimes_H \idd_P)(B)) =
(\sigma _{+n}\otimes_H \idd_P)(A*B),
                             \label{eq:3.4}
\end{equation}
where
$i\sigma _{+n}=i$
for
$i=1,\dots ,n$,
$(n+j)\sigma _{+n} =n+j\sigma $
for
$j=1,\dots ,m$.

A pseudoproduct
$*: P\otimes P\to (H\otimes H)\otimes _H P$
can be completely described by a family of
binary algebraic operations.
Let $P$ be an $H$-pseudoalgebra, $X=H^*$.
Lemma~\ref{lem 2.1} implies that for every
$a, b\in P$
their pseudoproduct has a unique presentation of the form
\[
a*b = \sum\limits_i (h_i\otimes 1)\otimes_{H} c_i,
\]
where $\{h_i\mid i\in I\}$ is a f\/ixed basis of $H$.
Consider the projections
(called Fourier coef\/f\/icients of~$a*b$)
\[
(a \oo{x} b) = \sum\limits_i \langle x, S(h_i) \rangle c_i \in P,
\]
for  all
$x\in X$.
The   $x$-products obtained have the following
properties:

\noindent {\bf  locality}
\begin{equation}
(a \oo{x_i} b) = 0 \quad \hbox{for almost all}\  i\in I;
                                              \label{eq3.5}
\end{equation}

\noindent {\bf sesqui-linearity}
\begin{equation}
(ha \oo{x} b) = (a \oo{xh} b),
\qquad
(a \oo{x} hb) = h_{(2)} (a_{S(h_{(1)})x} b).
                                              \label{eq3.6}
\end{equation}

Note that the locality property
does not depend on the choice of
a basis in $H$: \eqref{eq3.5}  means that
$\codim \{x\in X \mid (a\oo{x} b)=0\}<\infty $.

\begin{remark}[\cite{BDK}]\label{rem 3.1}
In the case $H=\Bbbk[D]$, $X\simeq \Bbbk[[t]]$,
where $\langle t^n, D^m\rangle  = n!\delta_{n,m}$,
the correspondence between conformal $n$-products ($n\ge 0$)
and the pseudoproduct is provided by
\[
a*b = \sum\limits_{n\ge 0}
\frac{1}{n!} ((-D)^{n} \otimes 1 )\otimes _H (a\oo{n} b),
\]
i.e.,
$a\oo{n} b = a\oo{t^n} b$, $n\ge 0$.
\end{remark}

In the same way, one may def\/ine Fourier coef\/f\/icients
of an arbitrary element $A\in H^{\otimes n}\otimes _H P$,
$n\ge 2$. By Lemma \ref{lem 2.1} $A$ can be uniquely
presented as
$A=\sum\limits_{\bar\imath}
(h_{i_1}\otimes \dots \otimes h_{i_{n-1}}\otimes 1)\otimes_H a_{\bar\imath}$,
$\bar\imath = (i_1,\dots, i_{n-1})\in I^{n-1}$.
By abuse of terminology, we will call these $a_{\bar\imath}\in P$
Fourier coef\/f\/icients of~$A$.

There is a canonical way to associate an ordinary algebra
$\mathcal A(P)$
with a given pseudoalgeb\-ra~$P$~\cite{BDK}.
As a linear space, $\mathcal A(P)$ coincides with
$X\otimes _H P$, and the product is given by
\[
 (x\otimes _H a)(y\otimes _H b ) = S^*(x_{(1)})y\otimes _H (a\oo{x_{(2)}} b),
 \qquad x,y\in X,\qquad a,b\in P.
\]
The algebra $\mathcal A(P)$ obtained is called the annihilation algebra
of~$P$.  If $P$ is a torsion-free $H$-module then the structure of $P$ can
be reconstructed from $\mathcal A(P)$~\cite{BDK}.

In the case of conformal algebras ($H=\Bbbk[D]$), there is a slightly dif\/ferent
construction called coef\/f\/icient algebra \cite{K1, K2, Ro1}.
Suppose $C$ is a conformal algebra and consider the space
$\Coeff C= \Bbbk[t,t^{-1}]\otimes _{\Bbbk[D]} C $,
where and $D$ acts on $\Bbbk[t,t^{-1}]$
as $t^nD = -nt^{n-1}$.
Denote $t^n\otimes _{\Bbbk[D]} a$ by
$a(n)$, $a\in C$, $n\in \mathbb Z$.
The product on $\Coeff C$ is provided by
\[
a(n)b(m) = \sum\limits_{s\ge 0} \binom{n}{s} (a\oo{s} b)(n+m-s),
\qquad n,m\in \mathbb Z, \qquad a,b\in C.
\]
An arbitrary conformal algebra
can be embedded into a conformal algebra of formal power series
over its coef\/f\/icient algebra \cite{K1}.

\subsection{Varieties of pseudoalgebras}\label{subsec3.2}

Suppose $\Omega $ is a variety of ordinary algebras over a f\/ield
of zero characteristic. Then $\Omega $ is def\/ined by a family of
homogeneous polylinear identities. Such an identity can be
written as
\begin{equation}
\sum\limits_{\sigma \in S_n}
t_\sigma (x _{1\sigma},\dots ,x _{n\sigma}) = 0,
                                                  \label{eq3.11}
\end{equation}
where each
 $t_\sigma(y_1,\dots ,y_n)$
is a linear combination
of non-associative words obtained from $y_1\dots y_n$
by some bracketing.

\begin{definition}[\cite{Kol}]\label{defn 3.3}
Let $\Omega $ be a variety of ordinary algebras
def\/ined by a family of homogeneous polylinear identities of the
form \eqref{eq3.11}.
Then set the $\Omega $ variety
 of pseudoalgebras as the class of pseudoalgebras
 satisfying the respective
``pseudo''-identities of the form
\begin{equation}
\sum\limits_{\sigma \in S_n}  (\sigma \otimes_H \idd_C)
t^*_\sigma (x_{1\sigma},\dots ,x_{n\sigma}) = 0,
                                                \label{eq3.12}
\end{equation}
where  $t_\sigma ^*$
means the same term $t_\sigma $ with respect to the
pseudoproduct operation~$*$.
\end{definition}

If $P$ is an $\Omega $ pseudoalgebra (or, in particular, conformal algebra)
then its annihilation (coef\/f\/icient) algebra belongs to
the $\Omega $ variety of ordinary algebras \cite{Kol}.
The converse is also true for conformal algebras.

However, the class of $\Omega $ pseudoalgebras is not a variety
in the ordinary sense: This class is not closed under
Cartesian products.

The main object of our study is the class of Jordan pseudoalgebras.
Recall that the variety of Jordan algebras is
def\/ined by the following identities:
\begin{equation}
a b=b a,
\qquad
((a a) b) a = (a a)(b a).
                                \label{eq3.13}
\end{equation}
In the polylinear form
(if $\hbox{char}\,\Bbbk \neq 2, 3$)
the second identity of
\eqref{eq3.13}
can be rewritten as follows
(see, e.g.,
\cite{ZSSS}):
\[
[abcd]  +  [dbca]
  +  [cbad]
= \{abcd\} + \{acbd\} + \{adcb\}.
\]
Here
$[\dots ]$
and
$\{\dots \}$
stand for the following bracketing schemes:
$[a_1\dots a_n]=(a_1[a_2\dots a_n])$,
$\{a_1a_2a_3a_4\}=((a_1a_2)(a_3a_4))$.

Therefore, an $H$-module $P$ (over a cocommutative Hopf algebra $H$)
endowed with a pseudoproduct~$\circ$
is a Jordan pseudoalgebra
if\/f it satisf\/ies the following identities of the form \eqref{eq3.12}:
\begin{gather}
a\circ b  =  (\sigma _{12}\otimes_H \idd_P)(b\circ a),
\nonumber\\
{[a\circ b\circ c\circ d] }
     +
 (\sigma _{14}\otimes _H \idd_P)[d\circ b\circ c\circ a]
  +
  (\sigma _{13}\otimes _H \idd_P)[c\circ b\circ a\circ d]
                      \nonumber\\
\qquad{}  = \{a\circ b\circ c\circ d \}
 + (\sigma _{23}\otimes _H \idd_P)\{a\circ c\circ b\circ d \}
+  (\sigma _{24}\otimes _H \idd_P)\{a\circ d\circ c\circ b \},
                                \label{eq3.16}
\end{gather}
where
$\sigma _{ij} = (i\,j)$
are  the transpositions from~$S_4$.

As in the case of ordinary algebras, the natural relations hold between
associative, Lie, and Jordan pseudoalgebras.
An associative pseudoalgebra $P$ with respect
to the new pseudoproduct
\[
[a* b]
= a*b - (\sigma _{12}\otimes_H \idd_P)(b*a)
\]
is a Lie pseudoalgebra denoted by $P^{(-)}$ \cite{BDK}.
Similarly,  another pseudoproduct $\circ $ given by
\[
a\circ b = a*b + (\sigma _{12}\otimes_H \idd_P)(b*a)
\]
makes $P$ into a Jordan pseudoalgebra $P^{(+)}$ \cite{Kol}.

\begin{example}\label{exmp 3.1}
Let $H'$ be a Hopf subalgebra of $H$,
and let $P'$ be an $H'$-pseudoalgebra
with respect to a pseudoproduct~$*'$.
Def\/ine a pseudoproduct on
$P=H\otimes _{H'} P'$
by linearity:
\[
(h\otimes_{H'} a)*(g\otimes_{H'} b) = \sum\limits_i (hh_i\otimes
gg_i)\otimes_{H} (1\otimes_{H'}c_i), \qquad g,h\in H,
\]
where
$a*'b = \sum\limits_i (h_i\otimes g_i) \otimes_{H'} c_i$,
$a,b,c_i \in P'$.
The pseudoalgebra $P$ obtained is called the current pseudoalgebra
$\Curr_{H'}^{H} P'$.

In particular,
 $\Bbbk \subset H$
is a Hopf subalgebra of~$H$. Hence, an ordinary algebra $A$
gives rise to current pseudoalgebra $\Curr A = \Curr_\Bbbk ^H A$.
\end{example}

It is clear that if $P'$ is an $\Omega $ pseudoalgebra over $H'$ then
so is $\Curr_{H'}^H P'$.

\begin{example}\label{exmp 3.2}
Consider $H=U(\mathfrak h)$,
where $\mathfrak h$ is a Lie algebra.
Then the free left $H$-module
$H\otimes H$ equipped by pseudoproduct
\[
(h\otimes a)*(g\otimes b)
 = (hb_{(1)}\otimes g)\otimes _H (1\otimes ab_{(2)}),
\qquad  a,b,g,h\in H,
\]
is an associative pseudoalgebra.
The submodule $W(\mathfrak h)=H\otimes \mathfrak h$ is a subalgebra of the
corresponding Lie pseudoalgebra $(H\otimes H)^{(-)}$.

Note that if $\mathfrak h' $ is a Lie subalgebra of $\mathfrak h$, then $H'=U(\mathfrak h')$
is a Hopf subalgebra of~$H$, and $\Curr_{H'}^H W(\mathfrak h')$ is actually a
subalgebra of $W(\mathfrak h)$.
\end{example}

In particular, if $\mathfrak h$ is the 1-dimensional Lie algebra
then $W(\mathfrak h)$ is just the Virasoro conformal algebra~\cite{K1}.

Later we will use the
classif\/ication of simple f\/inite Lie
pseudoalgebras~\cite{BDK}.
Although the results obtained in
\cite{BDK} are much more explicit, the
following statements are suf\/f\/icient for our purposes.

\begin{theorem}[\cite{BDK}]\label{thm 3.2}
A simple
finite
Lie pseudoalgebra  $L$
over
$H=U(\mathfrak h)$, $\dim \mathfrak h<\infty $,
$\Bbbk= \mathbb C$,
is isomorphic either to
 $\Curr \mathfrak g$, where
 $\mathfrak g$ is a simple finite-dimensional Lie algebra,
 or to a subalgebra of
$W(\mathfrak h)$.
\qed
\end{theorem}

\begin{theorem}[\cite{BDK}]\label{thm 3.3}
A simple
 Lie
 pseudoalgebra $L$
 over $H=U(\mathfrak h)\# \Bbbk[\Gamma]$
which is finite over $U(\mathfrak h)$
($\dim \mathfrak h<\infty $, $\Bbbk=\mathbb C$)
is a finite direct sum of isomorphic  simple
$U(\mathfrak h)$-pseudoalgebras
such that $\Gamma $ acts on them transitively.
\qed
\end{theorem}

\subsection{Conformal linear maps}\label{subsec3.3}

Let  $H$ be a cocommutative Hopf algebra, and let
 $M_1$, $M_2$ be two left $H$-modules.
A map
 $\varphi  :  M_1\to (H\otimes H )\otimes_H M_2$
is said to be (left)
conformal linear if
\[
\varphi (ha)=(1\otimes h)\varphi (a), \qquad h\in H,\qquad
a\in M_1.
\]
The set of all left conformal linear maps is
denoted by
$\Chom^l(M_1, M_2)$. For $M_1=M_2=M$ we denote
 $\Chom^l(M,M) = \Cend^l(M)$.

For every $H$-modules $M_1$, $M_2$,
the set $\Chom^l(M_1,M_2)$ can be considered
as an $H$-module with respect to the action
\[
h\varphi(a) = (h\otimes 1)\varphi(a),
\qquad
h\in H,\qquad \varphi \in \Chom^l(M_1,M_2),\qquad
a\in M_1.
\]

For example, if $P$ is a pseudoalgebra, $a\in P$,
then the operator of left multiplication
$L_a :  b\mapsto a*b$, $b\in P$,
belongs to $\Cend^l(P)$.

In order to unify notations, we will use
$\varphi * a$ for $\varphi (a)$,
$a\in M$, $\varphi \in \Cend^l(M)$.
One may consider
$*$ here as an $H$-bilinear map  from
$\Cend^l(M)\otimes M$ to $(H\otimes H)\otimes_H M$.
The relation \eqref{eq3.2} allows to expand this map to
\[
(H^{\otimes n}\otimes_H \Cend^l(M))\otimes
 (H^{\otimes m}\otimes_H M) \to H^{\otimes n+m}\otimes_H M.
\]
The correspondence between
$\varphi * a$ and $(\varphi \oo{x} a)$ ($x\in X$)
is given by
\[
  \varphi * a=
 \sum\limits_{i\in I} (S(h_i)\otimes 1)\otimes_H (\varphi \oo{x_i} a).
\]

The space $\Cend^l(M)$
can be also endowed with a family of $x$-products
given by
\[
(\varphi \oo{x} \psi)\oo{y} a = (\varphi \oo{x_{(2)}}
  (\psi \oo{S^*(x_{(1)})y} a)),
\qquad \varphi, \psi\in \Cend^l(M),\qquad
x,y\in X,
\]
for $a\in M$.

The $x$-products  $(\cdot \oo{x}\cdot )$ on  $\Cend^l(M)$
satisfy  \eqref{eq3.6},
but \eqref{eq3.5} does not hold, in general.
To ensure the locality, it is suf\/f\/icient
to assume that $M$ is a f\/initely generated $H$-module~\cite[Section~10]{BDK}.
Therefore, $\Cend^l(M)$ for a f\/initely generated
$H$-module $M$ is an associative $H$-pseudoalgebra.

For a f\/inite  pseudoalgebra $P$, it is easy to
rewrite the identity \eqref{eq3.16} using the
operators of left multiplication.
Namely, this identity is equivalent to
\begin{gather*}
  L_a*L_b*L_c + (\sigma _{13}\otimes _H \idd)
  (L_c*L_b*L_a) +
(\sigma _{123}\otimes _H \idd) L_{b*(c*a)}
\nonumber\\
\qquad{} = L_{a*b}*L_c + (\sigma _{23}\otimes _H \idd)
   (L_{a*c}*L_b)  +
  (\sigma _{13}\otimes _H \idd) (L_{c*b}*L_a),
\end{gather*}
where $\sigma _{123}$ denotes the
permutation $(1\,2\,3)\in S_4$.

\section[Tits-Kantor-Koecher construction for finite
Jordan pseudoalgebras]{Tits--Kantor--Koecher construction\\ for f\/inite
Jordan pseudoalgebras}\label{sec4}

The general scheme described in
\cite{Kantor, Koeher1, Koeher2, T}
provides an embedding of a Jordan algebra into a Lie algebra.
It is called the Tits--Kantor--Koecher (TKK) construction
for Jordan algebras.

Let us recall the TKK construction for ordinary algebras.
For a Jordan algebra $\mathfrak j$, the set
of derivations
$\Derr(\mathfrak j)$
 is a Lie subalgebra of
$\End(\mathfrak j)$
with respect to the commutator of linear maps.
Consider the (formal) direct sum
$\mathrm S(\mathfrak j)=\Derr(\mathfrak j)\oplus L(\mathfrak j)$,
where
$L(\mathfrak j)$
is the linear space of all left multiplications
$L_a: b\mapsto ab$, $a\in \mathfrak j$.
It is well-known that
$[L(\mathfrak j), L(\mathfrak j)] \subseteq \Derr(\mathfrak j)$.
Then the space~$\mathrm S(\mathfrak j)$
with respect to the new operation $[\cdot, \cdot]$
given by
\[
[(L_a+D),(L_b+T)]=L_{Db}-L_{Ta} + [L_a,L_b] + [D,T].
\]
is a Lie algebra called the structure Lie
algebra of~$\mathfrak j$.
Finally, consider
\[
\mathrm T (\mathfrak j)= {\mathfrak j}^- \oplus \mathrm S_0(\mathfrak j) \oplus
{\mathfrak j}^+,
\]
where
${\mathfrak j}^{\pm}\simeq \mathfrak j$,
$\mathrm S_0(\mathfrak j)$
is the subalgebra of
$\mathrm S(\mathfrak j)$,
generated by
$U_{a,b}  = L_{ab} +  [L_a, L_b]\in \mathrm S(\mathfrak j)$,
$a,b\in \mathfrak j$.
Let us endow
$\mathrm T(\mathfrak j)$ with the following operation:
\begin{gather*}
  [\Sigma ,a^-] = (\Sigma a)^-, \qquad
  [a^- , b^+] = U_{a,b},
                \\
      [a^+ , b^+] = [a^- , b^-] = 0,
\qquad
  [a^-, \Sigma ] = - (\Sigma  a)^-,   \qquad
  [a^+, \Sigma ]  = -(\Sigma ^* a)^+,
\\
 [\Sigma ,a^+]  = (\Sigma ^* a)^+,  \qquad
  [a^+ , b^-] =U^*_{a,b},
\end{gather*}
where
$\Sigma ^* = -L_a+ D$
for
 $\Sigma =L_a + D \in \mathrm S(\mathfrak j)$.
This operation makes $\mathrm T(\mathfrak j)$
to be a Lie algebra
 called the TKK construction for~$\mathfrak j$.

In the case of conformal algebras,
a similar construction was introduced in~\cite{Z2000}
by making use of coef\/f\/icient algebras.
We are going to get an analogue of TKK construction
for f\/inite Jordan pseudoalgebras
using the language of pseudoalgebras rather than
annihilation algebras.

\begin{definition}\label{defn 4.1}
Let $P$ be an $H$-pseudoalgebra.
A conformal endomorphism $T\in \Cend^l(P)$
is said to be a (left) {\em pseudoderivation},
if
\begin{equation}
T*(a*b) = (T*a)*b
+ (\sigma_{12}\otimes_H \idd_P)(a*(T*b))
                                              \label{eq4.1}
\end{equation}
for all
$a,b\in P$.
The set of all pseudoderivations of $P$ we denote by
$\Derr^l(P)$.
\end{definition}

In particular, if $P$ is a f\/inite pseudoalgebra
then \eqref{eq4.1}
is equivalent to
$[T*L_a]=L_{T*a}$, $a\in P$.

\begin{lemma}\label{lem 4.1}
Suppose that for
some $A\in (H\otimes H)\otimes_H \Cend^l(P)$
the equality
\[
A*(a*b) = (A*a)*b + (\sigma_{132}\otimes_H\idd)(a*(A*b))
\]
holds for all $a,b\in P$.
Then all Fourier coefficients of $A$
 belong to in $\Derr^l(P)$.
\end{lemma}

\begin{proof}
For every
$B\in H^{\otimes n+1}\otimes_H M$
($M$ is an $H$-module),
there exists a unique presentation
\[
B=\sum\limits_{\bar \imath}
 (G_{\bar  \imath}\otimes 1)\otimes_H b_{\bar  \imath},
\]
where $G_{\bar  \imath}$
form a linear basis of $H^{\otimes n}$
(see Lemma \ref{lem 2.1}).
By $B_{x_1,\dots, x_{n}}$,
 $x_i\in X$, we denote the expression
\[
\sum\limits_{\bar \imath} \langle x_1\otimes \dots \otimes x_n,
                     G_{\bar  \imath}\rangle b_{\bar  \imath}.
\]
It is clear that the map
$(x_1,\dots, x_n)\mapsto B_{x_1,\dots, x_n}$
is polylinear.
If we f\/ix an arbitrary set of $n-2$ arguments,
then the map
$X\otimes X \to M$
obtained
is local in the
sense of Def\/inition~\ref{defn 2.2}.

Let
$A=\sum\limits_{i\in I}
 (h_i\otimes 1)\otimes _H D_i$,
$a*b =
 \sum\limits_{j\in I} (h_j\otimes 1)\otimes _H c_j$,
$D_i*c_j =
  \sum\limits_{k\in I} (h_k\otimes 1)\otimes _H d_{ijk}$.
Then
\begin{gather}
D_i*(a*b)
  = \sum\limits_{j,k\in I}
  (h_k\otimes h_j\otimes 1)\otimes _H d_{ijk},
                                    \label{eq4.2} \\
A*(a*b) = \sum\limits_{i,j,k\in I}
( h_ih_{k(1)}\otimes h_{k(2)}\otimes h_j
 \otimes 1)\otimes _H d_{ijk}.
                                         \label{eq4.3}
\end{gather}
Compare \eqref{eq4.2} and \eqref{eq4.3} to get
\begin{equation}
(A*(a*b))_{x, y, z}
 =\sum\limits_{i\in I}\langle x_{(1)},h_i\rangle
 (D_i*(a*b))_{x_{(2)}y,z}.
                                     \label{eq4.4}
\end{equation}
In the same way,
\begin{gather}
 ((A*a)*b)_{x, y, z}
 =  \sum\limits_{i\in I}\langle x_{(1)},h_i\rangle
 ((D_i*a)*b)_{x_{(2)}y,z},
                                     \label{eq4.5}\\
 \big( (\sigma_{132}\otimes_H\idd)
 (a*(A*b))\big)_{x, y, z}
  =   \sum\limits_{i\in I}\langle x_{(1)},h_i\rangle
  \big( (\sigma_{12}\otimes_H\idd)
 (a*(D_i*b))\big)_{x_{(2)}y,z}.
                                     \label{eq4.6}
\end{gather}

The relations \eqref{eq4.4}--\eqref{eq4.6}
 together with Lemma~\ref{lem 2.3} imply
\[
\pi(x,y,z) =  \sum\limits_{i\in I}
\langle x,h_i\rangle
\big(
D_i*(a*b)
 -  (D_i*a)*b) -
(\sigma_{12}\otimes_H\idd)((D_i*a)*b)
\big )_{y,z} = 0.
\]
It means that
\begin{gather*}
D_i*(a*b)  - (D_i*a)*b) -
(\sigma_{12}\otimes_H\idd)(a*(D_i*b)) =0.\tag*{\qed}
\end{gather*}\renewcommand{\qed}{}
\end{proof}

\begin{lemma}\label{lem 4.2}
For  a finite pseudoalgebra $P$
the set of all pseudoderivations
is a subalgebra of  the Lie pseudoalgebra $\Cend^l(P)^{(-)}$.
\end{lemma}

\begin{proof}
Let $D_1, D_2\in \Derr^l(P)$, i.e.,
$[D_i*L_a]=L_{D_i*a}$ for $a\in P$,
$i=1,2$.

Since $\Cend^l (P)^{(-)}$ satisf\/ies Jacobi identity,
\begin{gather*}
[[D_1*D_2]*L_a]
 =
   [D_1*[D_2*L_a]]
   -(\sigma_{12}\otimes_H\idd)([D_2*[D_1*L_a]])  \\
\phantom{[[D_1*D_2]*L_a]}{}  =  [D_1*L_{D_2*a}]
  - (\sigma_{12}\otimes_H\idd)([D_2*L_{D_1*a}]) \\
\phantom{[[D_1*D_2]*L_a]}{}  =   L_{D_1*(D_2*a)}
   - (\sigma_{12}\otimes_H\idd)L_{D_2*(D_1*a)}=
   L_{[D_1*D_2]*a}.
             \nonumber
\end{gather*}
Hence, for every
$a,b\in P$ we have
\[
 [D_1*D_2]*(a*b) =  ([D_1*D_2]*a)*b
                    +
 (\sigma_{132}\otimes_H\idd)(a*([D_1*D_2]*b)).
\]
Lemma \ref{lem 4.1} implies that
$[D_1\oo{x} D_2]\in \Derr^l(P)$
for all $x\in X$.
\end{proof}

\begin{lemma}\label{lem 4.3}
Let $J$ be a finite Jordan pseudoalgebra,
and let
  $L(J)$ be the $H$-submodule of $\Cend^l(J)$
generated by
 $\{L_a\mid a\in J\}$.
Then
$[L_a \oo{x} L_b]$
is a pseudoderivation for every $x\in X$, i.e.,
$L'(J)\subseteq \Derr^l(J)$.
\end{lemma}

\begin{proof}
The following relation is easy to deduce from \eqref{eq3.16}:
\begin{gather*}
 L_a*L_{c*d}  +
  (\sigma _{13}\otimes _H \idd)
(L_d*L_{c*a}) +
(\sigma _{12}\otimes _H \idd) (L_c*L_{a*d})\\
 \qquad{} =  (\sigma _{123}\otimes _H \idd)(L_{c*d}*L_a) + L_{a*c}*L_d +
 (\sigma _{23}\otimes _H \idd) (L_{a*d}*L_c).
                                            \nonumber
\end{gather*}
So by  \eqref{eq:3.3}, \eqref{eq:3.4}
\begin{equation}
[L_a*L_{c*d}] = [L_{a*c}*L_d]
- (\sigma _{12}\otimes _H \idd)[L_c*L_{a*d}].
                                         \label{eq4.8}
\end{equation}

It is suf\/f\/icient to prove that for every
$a,b,c\in J$
we have
$[[L_a*L_b]*L_c] = L_{[L_a*L_b]*c}$,
i.e.,
\begin{gather}
  L_a*L_b*L_c   - (\sigma _{132}\otimes _H \idd)
  (L_c*L_a*L_b)
  +  (\sigma _{13}\otimes _H \idd)  L_c*L_b*L_a
\nonumber \\
 \qquad{} - (\sigma _{12}\otimes _H \idd)  L_b*L_a*L_c =L_{a*(b*c)} -
  (\sigma _{12}\otimes _H \idd)L_{b*(a*c)}.
                                                    \label{eq4.9}
\end{gather}
Indeed,
\begin{gather}
   L_a*L_b*L_c + (\sigma _{13}\otimes _H \idd)
   (L_c*L_b*L_a)  =
   - (\sigma _{123}\otimes _H \idd) L_{b*(c*a)}
 \nonumber\\
\qquad\quad{}   +   L_{a*b}*L_c + (\sigma _{23}\otimes _H \idd) (L_{a*c}*L_b)
   + (\sigma _{13}\otimes _H \idd) (L_{c*b}*L_a),
                                                    \label{eq4.10}
\\
  (\sigma _{12}\otimes _H \idd)  (L_b*L_a*L_c)
  +  (\sigma _{132}\otimes _H \idd)  (L_c*L_a*L_b)  \nonumber\\
\qquad{}  =(\sigma _{12}\otimes _H \idd)  (L_b*L_a*L_c)
  +(\sigma _{12}\sigma _{13}\otimes _H \idd)
    (L_c*L_a*L_b) \nonumber\\
\qquad{}  = - (\sigma _{12}\sigma _{123}\otimes _H \idd) L_{a*(c*b)}
 + (\sigma _{12}\otimes _H \idd) (L_{b*a}*L_c)  \nonumber\\
\qquad\quad{}  + (\sigma _{12}\sigma _{23}\otimes _H \idd) (L_{b*c}*L_a)
  +(\sigma _{12}\sigma _{13}\otimes _H \idd) (L_{c*a}*L_b)
  \nonumber \\
 \qquad{} = - (\sigma _{23}\otimes _H \idd) L_{a*(c*b)}
  +   (\sigma _{12}\otimes _H \idd) (L_{b*a}*L_c) \nonumber\\
  \qquad\quad{}+ (\sigma _{123}\otimes _H \idd) (L_{b*c}*L_a)  +
  (\sigma_{132}\otimes _H \idd) (L_{c*a}*L_b).
                                          \label{eq4.11}
\end{gather}

Subtracting
\eqref{eq4.10} from \eqref{eq4.11}
and using commutativity
$L_{a*b}=(\sigma _{12}\otimes _H \idd) L_{b*a}$,
we ob\-tain~\eqref{eq4.9}.
\end{proof}

\begin{definition}\label{defn 4.2}
Let $J$ be a f\/inite Jordan pseudoalgebra.
The formal direct sum of $H$-modules
\[
\mathrm S(J)= L(J) \oplus \Derr^l(J)
\]
endowed with the pseudoproduct
\begin{equation}
[(L_a+D)*(L_b+T)]=L_{D*b}-(\sigma _{12}\otimes_H \idd)L_{T*a} +
[L_a*L_b] + [D*T]
                                               \label{eq4.12}
\end{equation}
is called the {\em structure Lie pseudoalgebra\/}
of~$J$.
\end{definition}

It is straightforward to check that the (pseudo) anticommutativity
and Jacobi identities hold for~\eqref{eq4.12}.

Consider the elements
$U_{a,b}
 = L_{a*b} +  [L_a* L_b]
 \in (H\otimes H)\otimes_H \mathrm S(J)$,
$a,b\in J$.
By
 $U_{(a\oo{x}b)}=L_{(a\oo{x} b)} +
  [L_a\oo{x}L_b]$,
$x\in X$,
we denote the Fourier coef\/f\/icients of $U_{a,b}$.
The linear space
$\mathrm S_0(J)$
generated  by the set
$\{ U_{(a\oo{x}b)} \mid a,b\in J, x\in X \}$
is an
 $H$-submodule of
 $\mathrm S(J)$.

\begin{proposition}\label{prop 4.1}
The $H$-module $\mathrm S_0(J)$
 is closed under the pseudoproduct \eqref{eq4.12},
 i.e., $\mathrm S_0(J)$
 is a~Lie pseudoalgebra.
\end{proposition}

\begin{proof}
Let us calculate
 $[U_{a,b}*U_{c,d}]$,
 $a,b,c,d\in J$.
Denote
 $D=[L_a*L_b]$,
$A=a*b$.
 Then
$[U_{a,b}*L_{c*d}]
  = [L_A*L_{c*d}] + L_{(D*c)*d} +
  (\sigma_{132}\otimes _H \idd) L_{c*(D*d)}$,
$[U_{a,b}*[L_c*L_d]]
 =[L_A*[L_c*L_d]] + [L_{D*c}*L_d] +
 (\sigma _{132}\otimes _H \idd) [L_c *L_{D*d}]$.
Therefore,
\[
[U_{a,b}*U_{c,d}]
 = [L_A*L_{c*d}] + [L_A*[L_c*L_d]] + U_{D*c,d}
 + (\sigma _{132}\otimes _H \idd) U_{c,D*d}.
\]
From the f\/irst summand of the right-hand side we obtain
$[L_A*L_{c*d}] = [L_{A*c}*L_d]
 - (\sigma _{132}\otimes _H \idd)[L_c*L_{A*d}]$
by using~\eqref{eq4.8}.
Moreover,
 $[L_A*[L_c*L_d]]
  = L_{(A*c)*d} - (\sigma _{132}\otimes_H \idd)
  L_{c*(A*d)}$.
Hence,
\begin{gather*}
 [U_{a,b}*U_{c,d}] =
 U_{A*c,d}+U_{D*c,d} +
  (\sigma _{132}\otimes_H \idd)
(U_{c,D*d} - U_{c,A*d}).\tag*{\qed}
\end{gather*}\renewcommand{\qed}{}
\end{proof}

Denote
$U^*_{a,b}=-L_{a*b} + [L_a*L_b]$,
$a,b\in J$.
Note that
$U_{a,b}^*
 = -(\sigma_{12} \otimes_H \idd)U_{b,a}$,
so all Fourier coef\/f\/icients of
$U^*_{a,b}$
lie in $\mathrm S_0(J)$.
If $J$ is a Jordan pseudoalgebra and $J^2 = J$, i.e.,
every
$a\in J$
lies in the subspace generated by the set
$\{(b\oo{x}c) \mid b,c\in J, x\in X \}$,
then
  $\mathrm S_0(J) \supset L(J)$.

Indeed, for every
 $a,b\in J$ we have
 $U_{a,b}+(\sigma _{12}\otimes_H \idd_J) U_{b,a}
  = 2L_{a*b}$,
  so
 $L_{(a\oo{x}b)}\in \mathrm S_0(J)$.
 Hence, $L(J)=L(J^2) \subset \mathrm S_0(J)$.

Let us consider the direct sum of
$H$-modules
\[
\mathrm T(J) = J^- \oplus \mathrm S_0(J) \oplus J^+,
\]
where
$J^+$ and $J^-$
are isomorphic copies of~$J$.
Given
 $a\in J$
 (or $A\in H^{\otimes n}\otimes_H J$),
 we will denote by
 $a^\pm$
 (or $A^\pm$)
 the image of this element in~$J^\pm$
 (or $H^{\otimes n}\otimes_H J^\pm$).
Def\/ine a pseudoproduct on
$\mathrm T(J)$
by the following rule:
for
$a^\pm,b^\pm \in J^\pm$,
$\Sigma \in \mathrm S_0(J)$
 set
\begin{gather}
{} [a^+ * b^-] = U^*_{a,b}, \qquad
[a^- * b^+] = U_{a,b}, \qquad
[a^+ * b^+] = [a^- * b^-] = 0,  \nonumber\\
[a^-* \Sigma ] = - (\sigma _{12}\otimes_H \idd)(\Sigma *a)^-,
\qquad
     [\Sigma *a^-] = (\Sigma *a)^-,
\label{eq4.14}\\
  [a^+* \Sigma ]  = -(\sigma _{12}\otimes_H \idd) (\Sigma ^* * a)^+,
  \qquad  [\Sigma * a^+]  = (\Sigma ^* * a)^+ .\nonumber
\end{gather}
Set the pseudoproduct on $\mathrm S_0(J)$
to be the same as \eqref{eq4.12}.
Here we have used
$\Sigma ^* = -L_a+ D$
for
 $\Sigma =L_a + D \in \mathrm S(J)$.

Denote the projections of $\mathrm T(J)$
on
$J^+$, $J^-$, $\mathrm S_0(J)$
by
 $\pi_+$, $\pi_-$, $\pi_0$, respectively.
It is straightforward to check that
 $\mathrm T(J)$
is a Lie pseudoalgebra.
This is an analogue of the Tits--Kantor--Koecher
construction for an ordinary Jordan algebra.

Note that the structure pseudoalgebra
is a formal direct sum of the
corresponding $H$-modules, so
the condition
\[
\Sigma *b = 0 \quad \mbox{for all}\quad b\in J
\]
does not imply
 $\Sigma =0$ in $\mathrm S(J)$.
However,
if $\Sigma =L_a+D\in \mathrm S(J)$ and
 $[\Sigma *b^-]=[\Sigma * b^+] = 0$
in $\mathrm T(J)$
 for all
 $b\in J$, then
$a*b+D*b=0$ and $-a*b +D*b=0$ by \eqref{eq4.14}.
Therefore, $a*b=D*b=0$ for all $b\in J$, i.e.,
$\Sigma =0$ in $\mathrm S(J)$.

\begin{proposition}\label{prop 4.2}
Let $J$
be a simple finite Jordan pseudoalgebra.
Then
 ${\mathcal L} = \mathrm T(J)$
 is a simple finite Lie pseudoalgebra.
\end{proposition}

\begin{proof}
Suppose that there exists a non-zero proper ideal
$I \lhd {\mathcal L}$.
Let
\[
J_\pm
= \{a\in J \mid a^\pm = \pi_\pm(b)
\ \mbox{for some}\ b\in I\}.
\]
Since $J^2=J$,
we have
  ${\mathcal L} \supset L(J)$.
 Hence,
$J_\pm \lhd J$.

Analogously,
$J^0_\pm = \{a\in J \mid a^\pm \in I\cap J^{\pm} \}$
are also some ideals in~$J$.

1) Consider the case $J_+ = J_- = 0$ (hence,
$J^0_+ = J^0_- = 0$).
Since  $I\neq 0$,
there exists
$\Sigma =L_b + D \in \mathrm S_0(J)\cap I$, $\Sigma \neq 0$.
But
$[\Sigma * J^\pm] \subseteq
 H^{\otimes 2}\otimes_H J^0_\pm = 0$,
 so $\Sigma =0$ as we have shown above, which is a contradiction.

2) Let  $J_+ = J$, $J^0_-=0$.
Then for each
 $a\in J$
there exists
 $a^+ + \Sigma  + d^- \in I$.
 Consider
\[
[[(a^+ + \Sigma + d^-)*b^-]*c^-]
  =  [(U^*_{a,b}+(\Sigma *b)^-)* c^-]
                            =  (U^*_{a,b}*c)^- \in
H^{\otimes 3}\otimes_H J^0_-=0.
\]
For every
$a,b,c\in J$
we have
\begin{equation}                \label{eq:4.15}
-L_{a*b}*c + [L_a*L_b]*c=0.
\end{equation}

If
$a*b = \sum\limits_i (h_i\otimes 1)\otimes _H (a\oo{x_i} b)$,
then
$b*a = \sum\limits_i (1\otimes h_i)\otimes _H (a\oo{x_i} b)$
by commutativity. Therefore,
$L_{a*b} = \sum\limits_i (h_i\otimes 1)\otimes_H L_{a\oo{x_i} b}
  = (\sigma_{12}\otimes _H \idd_{L(J)})L_{b*a}$.
By the def\/inition of commutator, $(\sigma_{12}\otimes _H \idd_J)([L_a*L_b]*c)
 = -[L_b*L_a]*c$.
Relation \eqref{eq:4.15} implies
 $-L_{b*a}*c + [L_b*L_a]*c=0$ by symmetry.
Hence,
$0=(\sigma_{12}\otimes_H \idd_J) ( -L_{b*a}*c + [L_b*L_a]*c )
= -L_{a*b}*c - [L_a*L_b]*c$. Compare the last relation with \eqref{eq:4.15}
to get
$L_{a*b}*c=0$ for all
$a,b,c\in J$. Then the condition
$J^2=J$ implies $L(J)=0$, which is a contradiction.

3) The case  $J_- = J$, $J^0_+=0$
is completely analogous.

Hereby, if either of the ideals
 $J^0_\pm $
 is zero, then
at least one of the ideals $J_\pm \lhd J$
has to be zero, which is impossible.
Hence,
 $J^0_+ = J^0_- = J$, i.e.,
 $I\supset J^+, J^-$.
Since the whole pseudoalgebra ${\mathcal L}$
 is generated by
 $J^+\cup J^-$,
we have $I={\mathcal L}$.
\end{proof}

\section{Structure of simple
  Jordan pseudoalgebras}\label{sec5}

We have shown (Proposition \ref{prop 4.2}),
that for a simple f\/inite Jordan pseudoalgebra $J$
its TKK construction ${\mathcal L} = \mathrm T(J)$
is a simple f\/inite Lie pseudoalgebra.
This allows to describe simple Jordan pseudoalgebras
using
the classif\/ication
of simple Lie pseudoalgebras~\cite{BDK}.

\subsection[The case $H=U(\mathfrak h)$]{The case $\boldsymbol{H=U(\mathfrak h)}$}\label{subsec5.1}

Throughout this subsection,
$H$ is the universal enveloping Hopf algebra of
a f\/inite-dimensional Lie algebra
$\mathfrak h$
over~$\mathbb C$.

\begin{proposition}\label{prop 5.1}
Let $J$ be a simple finite Jordan
 $H$-pseudoalgebra.
Then the
  TKK construc\-tion~ $\mathrm T(J)$
is isomorphic to
the current algebra
$\Curr \mathfrak g$
over a simple finite-dimensional Lie al\-gebra~$\mathfrak g$.
\end{proposition}

\begin{proof}
If $J$
is a simple f\/inite Jordan $H$-pseudoalgebra, then
$\mathrm T(J)$ is a simple f\/inite Lie pseudoalgebra.
Hence, either
$\mathrm T(J) = \Curr \mathfrak g$,
where
$\mathfrak g$ is a simple f\/inite-dimensional Lie algebra, or
$\mathrm T(J)$
is a subalgebra in
$W(\mathfrak h)$
(see Theorem \ref{thm 3.2} and Example~\ref{exmp 3.2}).
The second case could not be realized
since by \cite[Proposition 13.6]{BDK}
the pseudoalgebra
$W(\mathfrak h)$
does not contain abelian subalgebras.
This is not the case for $\mathrm T(J)$.
\end{proof}

It remains to show that if
$\mathrm T(J)
 = J^+ \oplus \mathrm S_0(J) \oplus J^-
 = \Curr \mathfrak g$
then
 $J$
is the current pseudoalgebra over
a simple f\/inite-dimensional Jordan algebra.

Suppose
$e_1,\dots, e_n$ is a basis of $\mathfrak h$. Then the set of monomials
\[
 e^{(\alpha )}=e_1^{(\alpha_1)}\dots e_n^{(\alpha_n)}, \qquad
 \alpha=(\alpha_1,\dots,\alpha_n)\in \mathbb Z^n,\qquad \alpha_i\ge 0,
\]
where $e_i^{(\alpha_i)}= \frac{1}{\alpha_i!}e^{\alpha_i}$,
is a basis of $H$.
In order to simplify notation, we assume
$e^{(\alpha )}=0$ whenever $\alpha $ contains a negative component.

Denote $|\alpha | = \alpha_1+\dots + \alpha_n$.
We will use the standard deg-lex order on
the set of monomials of the form $e^{(\alpha)}$:
$e^{(\alpha)} \le e^{(\beta)}$ if
and only if $\alpha \le \beta$, i.e.,
either $|\alpha| < |\beta| $
or $|\alpha | = |\beta |$
and $\alpha $ is lexicographically less than~$\beta $.

Suppose the multiplication rule in $H$ is given by
$e^{(\alpha )}e^{(\beta )} =
\sum\limits_{\mu }
 \gamma _\mu ^{\alpha ,\beta } e^{(\mu )}$.
It is also useful to set $\gamma _\mu ^{\alpha ,\beta }=0$
if either of $\alpha$, $\beta$, $\mu$ contains a negative component.

The standard coproduct on $H$ is easy to compute in this notation:
$\Delta (e^{(\alpha)}) = \sum\limits_{\nu }e^{(\alpha-\nu)}\otimes e^{(\nu)}$.

\begin{theorem}\label{thm 5.1}
Let
$J$ be a simple finite
Jordan pseudoalgebra over
$H=U(\mathfrak h)$, where
$\mathfrak h$
is a finite-dimensional Lie algebra over the field $\mathbb C$.
Then
$J$
is isomorphic to the current algebra
$\Curr \mathfrak j$
over a finite-dimensional simple Jordan algebra~$\mathfrak j$.
\end{theorem}

\begin{lemma}\label{lem 5.1}
Let
$C=\Curr \mathfrak g = H\otimes \mathfrak g$.
Consider an arbitrary pair of elements
$a,b\in C$,
$a=\sum\limits_{\alpha } e^{(\alpha )}\otimes a_\alpha $,
$b=\sum\limits_{\beta } e^{(\beta )}\otimes b_\beta $,
$a_\alpha , b_\beta \in \mathfrak g$.
If
$[a*b]=0$, then
$[a_\alpha b_\beta ]=0$ in~$\mathfrak g$
for all
$\alpha$, $\beta $.
\end{lemma}

\begin{proof}
Straightforward computations show that
\begin{equation}
[a*b]
=\sum\limits_{\alpha ,\beta ,\nu ,\mu } (-1)^{|\beta -\nu |}
 \gamma _\mu ^{\alpha , \beta -\nu } \big(e^{(\mu )} \otimes 1\big)
 \otimes_H \big(e^{(\nu )}\otimes [a_\alpha b_\beta ]\big).
                                 \label{eq5.2}
\end{equation}
If
$[a*b]=0$
then \eqref{eq5.2}
implies that for every
$\nu ,\mu $
we have
\begin{equation}
\sum\limits_{\alpha ,\beta } (-1)^{|\beta -\nu |}
 \gamma _\mu^{\alpha , \beta -\nu }
 [a_\alpha b_\beta ] = 0.
                                          \label{eq5.3}
\end{equation}

Put
$\nu =\beta ^{\max}$
in \eqref{eq5.3}
(i.e.,
$b _\nu \ne 0$,
but
$b_\beta =0$ for all
$\beta >\nu $).
We obtain
$ \sum\limits_{\alpha }
 \gamma _\nu ^{\alpha ,0}
  [a_\alpha b_\nu ] = 0$
for each~$\mu $.
However,
\[
 \gamma _\mu ^{\alpha ,0}=
\begin{cases}
 0, & \mu \ne \alpha, \\
 1, & \mu =\alpha ,
\end{cases}
\]
hence,
$[a_\alpha b_{\beta ^{\max}}]=0$ for each~$\alpha $.

To f\/inish the proof, use the induction on~$\beta $.
Suppose that
$[a_\mu b_\beta ]=0$
for all
$\mu $ and
\mbox{$\beta >\beta _0$}.
Let us show that $[a_\mu b_{\beta _0}]=0$.
Put $\nu =\beta _0$.
Relation \eqref{eq5.3} implies
$0= \sum\limits_{\alpha }
 \gamma _\mu ^{\alpha ,0} [a_\alpha b_{\beta _0}]
 + \sum\limits_{\beta >\beta _0} (-1)^{|\nu -\beta |}
  \gamma _\mu^{\alpha , \beta -\nu }
  [a_\alpha b_\beta ]$.
The second summand is equal to zero by the inductive assumption.
So we have
$[a_\alpha b_{\beta _0}]=0$ for each~$\alpha $.
\end{proof}

Now, let
${\mathcal L} = J^+ \oplus \mathrm S_0(J) \oplus J^-
 = \Curr \mathfrak g$.
By
$\mathfrak j_0^\pm$
we denote the spaces spanned by all coef\/f\/icients
 $a_\alpha \in \mathfrak g$ ingoing in
the sums
$\sum\limits_{\alpha } e^{(\alpha)} \otimes a_\alpha \in J^\pm$.
Lemma~\ref{lem 5.1}
implies the spaces
$\mathfrak j_0^\pm$
are Abelian subalgebras of~$\mathfrak g$ such that $[H\otimes \mathfrak j_0^\pm * J^\pm] = 0$.
Moreover,
$H\otimes \mathfrak j_0^\pm \supseteq J^\pm$.

\begin{lemma}\label{lem 5.2}
Let
${\mathcal L}=\mathrm T(J)=\Curr \mathfrak g$,
where $\mathfrak g$ is a finite-dimensional
Lie algebra. Suppose that
there are no non-zero ideals $I \lhd {\mathcal L}$
such that
$\pi^\pm(I)=0$.
Then $J^\pm  = H\otimes \mathfrak j_0^\pm$, respectively.
\end{lemma}

\begin{proof}
It is enough to consider the ``$+$'' case.
Consider an arbitrary element
$a\in H\otimes \mathfrak j_0^+$,
$a=\pi_+(a) + \pi_0(a) + \pi_-(a)$.
Denote
$J_0^- = \pi_- (H\otimes \mathfrak j_0^+)$,
$J_0^0 = \pi_0 (H\otimes \mathfrak j_0^+)$.

For every $b\in J^+$ we have
$0=[a*b]=[\pi_+(a)*b] + [\pi_0(a)*b] + [\pi_-(a)*b]$.
Since
$[\pi_0(a)*b]\in H^{\otimes 2}\otimes_H \mathrm S_0(J)$,
$[\pi_-(a)*b] \in H^{\otimes 2}\otimes_H J^+$,
$[\pi_+(a)*b]=0$,
then
\begin{equation}
[J_0^- * J^+] = [J_0^0 * J^+] = 0.
                                       \label{eq5.4}
\end{equation}

Given
$H$-submodules
$A,B \subseteq {\mathcal L}$,
denote by
 $[A\cdot B]\subseteq {\mathcal L}$
the $H$-module
spanned (over $\mathbb C$) by all Fourier coef\/f\/icients
of all elements from
$[A*B]$.
By $[A^\omega\cdot  B]$
we denote the sum of $H$-modules
$\sum\limits_{n\ge 0} [A^n \cdot B]$,
where
$[A^0\cdot B] = B$,
$[A^{n+1}\cdot B]=[A\cdot [A^n\cdot B]]$.

For example,
$[\mathrm S_0(J)^\omega \cdot J_0^-]\subseteq J^-$.
Moreover,
the Jacobi identity and \eqref{eq5.4}
imply
$[J^+ * [\mathrm S_0(J)^\omega\cdot  J_0^-]] =0$.
It is also easy to note that
$[\mathrm S_0(J)
  * [\mathrm S_0(J)^\omega\cdot  J_0^-]]
   \subseteq
  H^{\otimes 2}\otimes _H
   [\mathrm S_0(J)^\omega\cdot  J_0^-]$.
Since
$[J^- * [\mathrm S_0(J)^\omega \cdot J_0^-]]=0$,
then
$I=[\mathrm S_0(J)^\omega\cdot  J_0^-]$
is a proper ideal of~$\mathcal L$,
$I\supseteq J_0^-$
and
$\pi^+(I)=0$.
Hence,
$I=0$, and
$J_0^-=0$.

Further, let us consider
\begin{equation}
I=[\mathrm S_0(J)^\omega\cdot  J_0^0]
+ [\mathrm S_0(J)^\omega\cdot  [J^- \cdot J_0^0]]
\subseteq
\mathrm S_0(J)\oplus J^-.
                                          \label{eq5.5}
\end{equation}
It follows from  \eqref{eq5.4} that
$[J^+ * [\mathrm S_0(J)^\omega \cdot J_0^0]]=0$.
Moreover,
$[J^+ *
 [\mathrm S_0(J)^\omega\cdot  [J^- \cdot J_0^0]]]
 \subseteq
  H^{\otimes 2} \otimes _H
   [\mathrm S_0(J)^\omega \cdot J_0^0]]$.
Therefore,
$[J^+\cdot  I] \subseteq I$.
Since
$[\mathrm S_0(J)\cdot  I] \subseteq I$
by  construction,
and $[J^- \cdot I] \subseteq I$ by the Jacobi identity,
the ideal \eqref{eq5.5}
is proper in $\mathcal L$, so $J_0^0=0$.

We have proved that
$ \pi_-(H\otimes \mathfrak j^+_0)
 = \pi_0(H\otimes \mathfrak j^+_0) =0$.
Thus,
 $J^+ = H\otimes \mathfrak j^+_0$.
\end{proof}

Hence,
under the conditions of Lemma~\ref{lem 5.2}
one has
$J= H\otimes \mathfrak j$,
$\mathfrak j\simeq \mathfrak j_0^\pm$.
Now it is necessary to show that the
Jordan pseudoproduct on~$J$
may be restricted to an ordinary Jordan product
on~$\mathfrak j$.

\begin{proposition}\label{prop 5.2}
Let $J$ be a finite Jordan
$H$-pseudoalgebra
such that
$\Ann_l(J):= \{ a\in J \mid a*J=0\} = 0$.
Assume
that $J=H\otimes \mathfrak j$, where
$\mathfrak j$ is a linear space.
If ${\mathcal L }=\mathrm T(J)=\Curr \mathfrak g$
 then
$\mathfrak j$ has a~structure of ordinary Jordan algebra
such that
$\mathfrak g\simeq \mathrm T(\mathfrak j)$.
\end{proposition}

\begin{proof}
Let $a,b\in J$ be some elements of the form
$a=1\otimes \alpha $,
$b=1\otimes \beta $,
$\alpha ,\beta \in \mathfrak j$.
Then
$2 L_{a*b} = [a^-*b^+]
 + (\sigma _{12}\otimes_H \idd) [b^-*a^+]
 =(1\otimes 1 )\otimes_H (1\otimes [\alpha ^-\beta ^+] +
 1\otimes [\beta ^- \alpha ^+])$
(here
$\alpha ^\pm$
denote the images of
$\alpha \in \mathfrak j$
in
$\mathfrak j_0^\pm$).

Thus,
 $L_{a*b}= (1\otimes 1)\otimes_H (1\otimes s(\alpha ,\beta ))\in \mathcal L$,
 where
 $s(\alpha ,\beta )\in [\mathfrak j^- \mathfrak j^+]\subseteq \mathfrak g$.
Therefore,
 $L_{(a\oo{t^\nu } b)} = 0$
 for
$\nu=(\nu _1,\dots ,\nu _n) >(0,\dots, 0)$.
Here we have used the notation
$t^\nu =t_1^{\nu _1}\dots t_n^{\nu _n}$
for basic functionals in
$X=H^*$.

Since $L_x=0$ implies
$x=0$,
we have
\begin{equation}
 a*b = (1\otimes 1)\otimes_H c, \qquad c\in J.
                                             \label{eq5.6}
\end{equation}
Suppose that
$c=\sum\limits_{\mu } e^{(\mu )}\otimes \gamma_\mu $,
$\gamma_\mu \in \mathfrak j$.
Assume that
 the maximal
  $\mu = \mu_{\max}$
 such that
 $\gamma _\mu \ne 0$
is a multi-index greater than
$(0,\dots ,0)$. Then
\begin{equation}
 [(1\otimes s(\alpha ,\beta ))
  \oo{t^{\mu_{\max}}} (1\otimes \delta^- )]=0
                                              \label{eq5.7}
\end{equation}
for all
$\delta \in \mathfrak j$.
On the other hand,
$ [(1\otimes s(\alpha ,\beta )) \oo{t^{\mu_{\max}}}
 (1\otimes \delta^- )]
 = (c \oo{t^{\mu _{\max}}} (1\otimes \delta ))^- $.
It is easy to see that
 the relations \eqref{eq5.6}, \eqref{eq5.7}
 and the axioms of a pseudoalgebra imply
$(c \oo{t^{\mu _{\max}}} (1\otimes \delta ))
 = ((1\otimes \gamma _{\mu_{\max}})
 \oo{\varepsilon } (1\otimes \delta )) = 0 $,
i.e.,
$L_{1\otimes \gamma _{\mu _{\max}}}=0$.
Thus,
 $\gamma _{\mu _{\max}} = 0$, which is a contradiction.

We have proved that
$ (1\otimes \alpha )*(1\otimes \beta )
 = (1\otimes 1)\otimes_H (1\otimes \gamma (\alpha ,\beta ))$,
$\gamma (\alpha ,\beta ) \in \mathfrak j$.
This relation leads to an ordinary product on
$\mathfrak j$ def\/ined by the rule
$\alpha  \cdot \beta  = \gamma (\alpha ,\beta )$.
Then the pseudoalgebra
 $J$
is a~current pseudoalgebra over
$\mathfrak j$,
and $(\mathfrak j, \cdot)$ is necessarily a
simple f\/inite-dimensional Jordan algebra.
To complete the proof, it is enough to note that
$\mathrm T(\Curr \mathfrak j)\simeq \Curr \mathrm T(\mathfrak j)$.
For f\/inite-dimensional Lie algebras
$\mathfrak g_1$, $\mathfrak g_2$
the condition
$\Curr \mathfrak g_1\simeq \Curr \mathfrak g_2$
implies
$\mathfrak g_1\simeq \mathfrak g_2$.
\end{proof}

\begin{proof}[Proof of Theorem \ref{thm 5.1}]
Let $J$ be a simple f\/inite Jordan pseudoalgebra.
 Proposition \ref{prop 5.1} implies that
${\mathcal L}=\mathrm T(J)\simeq \Curr \mathfrak g$,
where
$\mathfrak g$ is a simple f\/inite-dimensional Lie algebra.
By Lemma \ref{lem 5.2},
$J=H\otimes \mathfrak j$.
Since $\mathcal L$ satisf\/ies the conditions
of Proposition \ref{prop 5.2},
we have $J\simeq \Curr\mathfrak j$,
$\mathrm T(\mathfrak j)=\mathfrak g$,
where
 $\mathfrak j $ is a simple f\/inite-dimensional
 Jordan algebra.
\end{proof}

\begin{corollary}[\cite{Z2000}]\label{cor 5.2}
A simple finite Jordan conformal algebra
is isomorphic to the current conformal algebra
over a simple finite-dimensional Jordan algebra.
\qed
\end{corollary}

\subsection[The case
  \protect{$H=U(\mathfrak h)\# \mathbb C[\Gamma ]$}]{The case
  $\boldsymbol{H=U(\mathfrak h)\# \mathbb C[\Gamma ]}$}\label{subsec5.2}

If $J$ is a pseudoalgebra over   $H=U(\mathfrak h)\# \mathbb C[\Gamma ]$
then it is in particular a pseudoalgebra over
$U(\mathfrak h)$. The structure of $H$-pseudoalgebra on $J$ is completely
encoded by $U(\mathfrak h)$-pseudoalgebra structure and by the action
of $\Gamma $ on $U(\mathfrak h)$, see \cite[Section~5]{BDK} for details.

\begin{theorem}\label{thm 5.2}
Let $J$ be a simple Jordan pseudoalgebra over
$H=U(\mathfrak h)\# \mathbb C [\Gamma ]$,
$\dim\mathfrak h<\infty$,
which is a finitely generated $U(\mathfrak h)$-module.
Then
\[
J\simeq \bigoplus\limits_{i=1}^{m} \Curr^{U(\mathfrak h)} \mathfrak j_i,
\]
where $\mathfrak j_i $
are isomorphic finite-dimensional
simple Jordan algebras, and
$\Gamma $
acts transitively on the family
$\{ \Curr^{U(\mathfrak h)} \mathfrak j_i: i=1,\dots, m \}$.
\end{theorem}

\begin{proof}
By Proposition \ref{prop 4.2}
${\mathcal L} = \mathrm T(J)$
is a simple $H$-pseudoalgebra, and it is clear that
$\mathcal L$ is a~f\/initely generated $U(\mathfrak h)$-module.
Theorem~\ref{thm 3.3} and Proposition~\ref{prop 5.1}
imply that
${\mathcal L}
  = \bigoplus\limits_{i=1}^{m} \Curr^{U(\mathfrak h)} \mathfrak g_i$
where
$\Curr^{U(\mathfrak h)} \mathfrak g_i = \Curr_i$
are isomorphic simple current Lie
$U(\mathfrak h)$-pseudoalgebras,
and  $\Gamma $ acts on them transitively.

Hence,
${ \mathcal L } = \Curr ^{U(\mathfrak h) }\tilde {\mathfrak g}$,
where
$\tilde {\mathfrak g} = \bigoplus\limits_{i=1}^m \mathfrak g_i$.
The $H$-pseudoalgebra $\mathcal L$ could
be considered as
an $U(\mathfrak h)$-pseudoalgebra endowed
with an action of  $\Gamma $ on it
which is compatible with that of $U(\mathfrak h)$:
$g(ha) = h^g (ga)$,
$h\in U(\mathfrak h)$, $a\in {\mathcal L}$,
$g\in \Gamma$.

Consider $\mathcal L$ as the current $U(\mathfrak h)$-pseudoalgebra
over $\tilde {\mathfrak g}$.
The condition of Lemma \ref{lem 5.2} holds for this $\mathcal L$.
Indeed,
if $I$ is an ideal of the $U(\mathfrak h)$-pseudoalgebra
$\mathcal L$ and  $\pi^\pm(I)=0$,
then $\Gamma I$ is a proper ideal
of $\mathcal L$ (as of an $H$-pseudoalgebra)
such that its projections onto $J^\pm$ are zero.
Moreover, if $J$ as an $H$-pseudoalgebra
has no non-trivial (left) annihilator
$\Ann_l(J)$  then
so is $J$ as an $U(\mathfrak h)$-pseudoalgebra
(see \cite[Corollary~5.1]{BDK}).

Therefore, the same arguments as in the proof
of Proposition \ref{prop 5.2} show that
$J = \Curr^{U(\mathfrak h)} \tilde{\mathfrak j}$,
where
$\tilde{\mathfrak j}$ is a f\/inite-dimensional Jordan algebra.

The explicit expression \cite[equation (5.7)]{BDK}
for pseudoproduct over $H$ shows that for every $x\in X=U(\mathfrak h)^*$,
$g\in \Gamma$, $a,b\in J$ we have
\[
(a\oo{x\otimes g^*} b)
 = (a \oo{(x\otimes 1)g} b) = (ga \oo{x} b),
\]
where
$\langle g^*, \gamma \rangle = \delta_{g,\gamma}$,
$\gamma \in \Gamma$.
Hence, the following relation between Fourier coef\/f\/icients of
$U_{a,b}$ holds:
$U_{(a\oo{x\otimes g^*}b)}= U_{(ga\oo{x} b)}$.
Here in the left- and right-hand sides we state Fourier coef\/f\/icients
over $H^*$ and $X$, respectively.
Therefore, the relations between the $H$-module
$\mathrm S_0({}_H J)$
 and the
$U(\mathfrak h)$-module ${}_{U(\mathfrak h)} \mathrm S_0(J)$
 are the same as between
$H$-module ${}_H J$ and $U(\mathfrak h)$-module
${}_{U(\mathfrak h)} J$.

Now it is clear that
$\tilde{\mathfrak g} = \mathrm T (\tilde{\mathfrak j})$.
Hence,
$\tilde{\mathfrak j} = \bigoplus\limits_{i=1}^m \mathfrak j_i$,
$\mathfrak g_i = \mathrm T (\mathfrak j_i)$,
where $\mathfrak j_i$
are simple Jordan algebras.
So,
$J = \bigoplus\limits_{i=1}^m \Curr^{U(\mathfrak h)} \mathfrak j_i$,
and $\Gamma $ necessarily acts on these
current algebras transitively.
\end{proof}

\subsection*{Acknowledgements}

I am very grateful to L.A.~Bokut, I.V.~L'vov,
E.I.~Zel'manov, and V.N.~Zhelyabin
for their interest in the present work
and helpful discussions.
It is my pleasure to appreciate the ef\/forts of the referees, whose suggestions
and comments helped me to make the paper readable.

The work was partially supported by SSc-344.2008.1.
I gratefully acknowledge the support
of the Pierre Deligne fund based on his 2004
Balzan prize in mathematics, and
Novosibirsk City Administration grant of 2008.

\pdfbookmark[1]{References}{ref}
\LastPageEnding


\begin{thebibliography}{99}

\footnotesize\itemsep=0pt

\bibitem{BDK}
Bakalov~B., D'Andrea~A., Kac~V.G.,
Theory of f\/inite pseudoalgebras, {\it Adv. Math.} {\bf 162} (2001), 1--140, \href{http://arxiv.org/abs/math.QA/0007121}{math.QA/0007121}.

\bibitem{BD}
Beilinson~A.A., Drinfeld~V.G.,
Chiral algebras,
{\it American Mathematical Society Colloquium Publications}, Vol.~51,
American Mathematical Society, Providence, RI, 2004.

\bibitem{BPZ}
Belavin~A.A., Polyakov~A.M., Zamolodchikov~A.B.,
Inf\/inite conformal symmetry in two-dimensional quantum f\/ield theory,
{\it Nuclear Phys.~B} {\bf 241} (1984), 333--380.

\bibitem{Bor}
Borcherds~R.E.,
Vertex algebras, Kac--Moody algebras, and the Monster,
{\it Proc. Nat. Acad. Sci. U.S.A.} {\bf 83} (1986), 3068--3071.

\bibitem{DK}
D'Andrea~A., Kac~V.G.,
Structure theory of f\/inite conformal algebras,
{\it Selecta Math. (N.S.)} {\bf 4} (1998), 377--418.

\bibitem{FK}
Fattori D., Kac V.G.,
Classif\/ication of f\/inite simple Lie conformal superalgebras,
{\it J. Algebra} {\bf 258} (2002),    23--59, \href{http://arxiv.org/abs/math-ph/0106002}{math-ph/0106002}.

\bibitem{K1}
Kac~V.G.,
Vertex algebras for beginners, 2nd ed.,
{\it University Lecture Series}, Vol.~10, American Mathematical Society, Providence, RI, 1998.

\bibitem{K2}
Kac~V.G.,
Formal distribution algebras and conformal algebras,
in Proceedinds of XII-th International Congress in Mathematical Physics (ICMP'97)
(Brisbane), Int. Press, Cambridge, MA, 1999, 80--97, \mbox{\href{http://arxiv.org/abs/q-alg/9709027}{q-alg/9709027}}.

\bibitem{KR}
Kac~V.G., Retakh~A.,
Simple Jordan conformal superalgebras,
{\it J. Algebra Appl.} {\bf 7} (2008),  517--533, \href{http://arxiv.org/abs/0801.0755}{arXiv:0801.0755}.

\bibitem{Kantor}
Kantor I.L.,
Classif\/ication of irreducible transitively dif\/ferential groups,
{\it  Soviet Math. Dokl.} {\bf 5} (1965), 1404--1407.

\bibitem{Koeher1}
Koecher~M.,
Embedding of Jordan algebras into Lie algebras.~I,
{\it Amer. J. Math.} {\bf 89} (1967), 787--816.

\bibitem{Koeher2}
Koecher~M.,
Embedding of Jordan algebras into Lie algebras.~II,
{\it Amer. J. Math.} {\bf 90} (1968), 476--510.

\bibitem{Kol}
Kolesnikov P.S.,
Identities of conformal algebras and pseudoalgebras,
{\it Comm. Algebra} {\bf 34} (2006), no.~6, 1965--1979, \href{http://arxiv.org/abs/math.RA/0412397}{math.RA/0412397}.

\bibitem{La}
Lambek~J.,
Deductive systems and categories~II.
 Standard constructions and closed categories,
 in Category Theory, Homology Theory and their Applications, I (Battelle Institute Conference, Seattle, Wash., 1968, Vol.~One),
Springer, Berlin, 1969, 76--122.

\bibitem{Ro1}
Roitman M.,
On free conformal and vertex algebras,
{\it J. Algebra} {\bf 217} (1999),   496--527, \href{http://arxiv.org/abs/math.QA/9809050}{math.QA/9809050}.


\bibitem{Sw2}
Sweedler M.E.,
Cocommutative Hopf algebras with antipode,
{\it Bull. Amer. Math. Soc.} {\bf 73} (1967),   126--128.

\bibitem{Sw}
Sweedler M.E.,
Hopf algebras,
{\it Mathematics Lecture Note Series},
W.A.~Benjamin, Inc., New York, 1969.

\bibitem{T}
Tits J.,
Une classe d'algebres de Lie en relation   avec algebres de Jordan,
{\it  Indag. Math.} {\bf 24} (1962), 530--535.

\bibitem{Z2000}
Zelmanov E.I.,
On the structure of conformal algebras,
in Proceedings of Intern. Conf. on Combinatorial and Computational Algebra,
(May 24--29, 1999, Hong Kong)
{\it Contemp. Math.} {\bf 264} (2000), 139--153.

\bibitem{ZSSS}
Zhevlakov K.A., Slin'ko A.M., Shestakov I.P., Shirshov A.I.,
Rings that are nearly associative,
{\it Pure and Applied Mathematics}, Vol.~104, Academic Press, Inc., New York~-- London, 1982.

\end{thebibliography}
\end{document}